\documentclass[a4paper,11pt]{article}

\usepackage{amsfonts}
\usepackage{amsmath}
\usepackage{amssymb}
\usepackage{graphicx}
\usepackage{amscd}
\usepackage{tikz}
\usepackage{epsfig}
\usepackage{psfrag}
\usepackage{geometry}
\usepackage{pgfplots}
\pgfplotsset{compat=1.15}
\usepackage{mathrsfs}
\usetikzlibrary{arrows}

\newtheorem{lemma}{Lemma}
\newtheorem{theorem}{Theorem}

\newtheorem{definition}{Definition}

\newtheorem{proposition}{Proposition}

\def\kraj{\hfill\rule{6pt}{6pt}}

\def\Z{\mathbb{Z}}

\def\dj{d\kern-0.4em\char"16\kern-0.1em}
\def \Dj {\mbox{\raise0.3ex\hbox{-}\kern-0.4em D}}
\definecolor{qqwuqq}{rgb}{0.,0.39215686274509803,0.}
\definecolor{wwccff}{rgb}{0.4,0.8,1.}
\definecolor{ffqqtt}{rgb}{1.,0.,0.2}
\definecolor{zzttqq}{rgb}{0.6,0.2,0.}
\definecolor{xdxdff}{rgb}{0.49019607843137253,0.49019607843137253,1.}
\definecolor{ffqqqq}{rgb}{1.,0.,0.}
\definecolor{zzttqq}{rgb}{0.6,0.2,0.}

\title{Lattice paths and quiver generating series with higher level generators}
\date{}
\author{Du\v san \Dj or\dj evi\'c\footnote{Faculty of Physics, University of Belgrade, Studentski Trg 16, 11000 Belgrade, Serbia.}
\and
Marko Sto\v si\'c\footnote{
CAMGSD, Departamento de Matem\'atica, 
Instituto Superior T\'ecnico,
Av. Rovisco Pais 1, 1049-001 Lisbon, Portugal, and
Mathematical Institute SANU, Knez Mihajlova 36, 11000 Belgrade, Serbia. 
(e-mail: {\tt mstosic@isr.ist.utl.pt}. Corresponding author.)
}}

\begin{document}
\maketitle
\begin{abstract}
The generalized knots-quivers correspondence extends the original knots-quivers correspondence, by allowing higher level generators of quiver generating series. In this paper we explore the underlined combinatorics of such generating series, relationship with the BPS numbers of a corresponding knot, and new combinatorial interpretations of the coefficients of generating series.
\end{abstract}
\section{Introduction}

Study of knots, using techniques from both mathematics and physics, led to some unexpected results. An example is the correspondence between knots and quivers (called knots-quivers correspondence, introduced in \cite{KRSS}), that associates to each knot a quiver so that the generating series of symmetrically colored HOMFLY-PT polynomials of a knot can be written as a quiver generating series. However, quiver generating series are series in (generically) more than one variable $x_i$ -- one for each vertex of the quiver -- and therefore differ from the generating series of HOMFLY-PT polynomials of a knot, that contain only one variable $x$. This implies that some relation between $x$ and $x_i$ has to be postulated. Starting from \cite{KRSSlong}, specialization of the form 
\begin{equation}\label{x}x_i=(-1)^{s_i}a^{a_i}q^{q_i}x,\quad i=1,\ldots,m,\end{equation}
was used, and it fruitfully reproduced many important results \cite{SW1,SW2}. However, there is no a priori reason why this relation has to be linear in $x$, and it was recently suggested in \cite{S} that one should generalize this relation in order to capture a more general class of knots. This generalization was also inspired by some other considerations (\cite{EGG+} contains a similar identification for $F_K$ invariants, see also \cite{EKL2} for another example), and it states that one should use identification of the form 
\begin{equation}\label{x1}x_i=(-1)^{s_i}a^{a_i}q^{q_i}x^{n_i},\quad i=1,\ldots,m,\end{equation}
where $n_i$ are non-negative integers, possibly greater than one.  This results in the so-called generalized knots-quivers correspondence, and the corresponding quivers where computed for many homologically thick knots  in \cite{S}.

In this paper we focus on the quiver generating series and the corresponding specializations where $n_i$ can be arbitrary non-negative integers, and we extend on the previous results where all specializations were linear. 

In particular, we shall be interested in the following problems: Let  $C$ be a symmetric quiver with $m$ vertices. To the $i$th vertex we associate a non-negative integer $n_i$. This $n_i$ we shall call a {\it{level}} of the vertex $i$.  We are interested in studying the limits of the quiver generating series
\begin{equation}
y(x_1,\ldots,x_m)=\lim_{q\to 1} \frac{P_C(q^{n_1}x_1,\ldots,q^{n_{j-1}}x_{j-1},q^{n_j} x_j,q^{n_{j+1}}x_{j+1},\ldots, q^{n_m}x_m)}{P_C(x_1,\ldots,x_m)}=\label{definv-bis0}
\end{equation}
$$=\sum_{l_1,\ldots,l_m}b_{l_1,\ldots,l_m}x_1^{l_1}\ldots x_m^{l_m},$$
in computing the coefficients $b_{l_1,\ldots,l_m}$, and in studying their combinatorial interpretation.

\section{Preliminaries}\label{prelim}

\subsection{BPS numbers}
For a knot $K$, let $P_r(a,q)$ denote the reduced $r$-th symmetrically colored HOMFLY-PT polynomial, and $\overline{P}_r(a,q)$ denote the unreduced $r$-th symmetrically colored HOMFLY-PT polynomial. Generating function of all colored {HOMFLY-PT} polynomials of a given knot $K$ can be re-expressed as the exponential in the following way

\begin{equation}\label{BPS}\overline{P}(x,a,q):=\sum_{r\ge 0} \overline{P}_r(a,q) x^r =\exp \left(  \sum_{n,r\ge 1} \frac{1}{n} f_r(a^n,q^n) x^{rn}  \right),
\end{equation}
where
$$f_r(a,q)=\sum_{i,j} \frac{N_{r,i,j}\,a^i q^j}{q-q^{-1}}.$$
The numbers $N_{r,i,j}$ are the BPS number. The famous LMOV conjecture \cite{LM,LMV,OV}, resulting from physical considerations, states that the numbers $N_{r,i,j}$ are integers.

Note that the formula (\ref{BPS}) can  be re-written equivalently as the multiple product:
\begin{equation*}
\overline{P}(x,a,q)=
\prod_{r> 0} \prod_{i,j\in\mathbb{Z}} \prod_{k\ge 0}\left(1 -   x^r  a^{\frac{i}{2}}q^{\frac{j+1}{2}+k} \right)^{N_{r,i,j}}.   
\end{equation*}

\subsection{Quivers and motivic Donaldson-Thomas invariants}
Let $Q$ be a quiver with $m$ vertices, and let $C$ be the corresponding adjacency matrix. Then the quiver generating series is given by 
\begin{equation}
P_C(x_1,\ldots,x_m)=\sum_{d_1,\ldots,d_m} \frac{(-q^{\frac{1}{2}})^{\sum_{i,j=1}^m C_{i,j}d_id_j}}{(q;q)_{d_1}\cdots(q;q)_{d_m}} x_1^{d_1}\cdots x_m^{d_m}. \label{P-C}
\end{equation}
The motivic Donaldson-Thomas invariants $\Omega_{d_1,\ldots,d_m;j}$ of a symmetric quiver $Q$  are encoded in the following product decomposition of the above series:
\begin{equation*}
P_C=
\prod_{(d_1,\ldots,d_m)\neq 0} \prod_{j\in\mathbb{Z}} \prod_{k\geq 0} \Big(1 -  \big( x_1^{d_1}\cdots x_m^{d_m} \big) q^{\frac{j+1}{2}+k} \Big)^{(-1)^{j+1}\Omega_{d_1,\ldots,d_m;j}}.   \label{PQx-Omega}
\end{equation*}

\begin{theorem}\cite{KS,Ef} The numbers $\Omega_{d_1,\ldots,d_m;j}$ are nonnegative integers. \end{theorem}
The numbers $\Omega_{d_1,\ldots,d_m;j}$ also have enumerative interpretation, like the intersection Betti numbers of the moduli space of all semisimple representations of $Q$, or as the Chow-Betti numbers of the moduli space of all simple representations \cite{MR,FR}.

\subsection{(Generalized) knots-quivers correspondence}
The knots-quivers correspondence was postulated in \cite{KRSS}. It states that for every knot $K$ there exists a symmetric quiver $Q$, such that the knot generating series $\overline{P}(x,a,q)$ matches the quiver generating series $P_C(x_1,\ldots,x_m)$, with the specialization $x_i=(-1)^{s_i} a^{a_i} q^{q_i} x$, for some integers $a_i$, $q_i$ and $s_i$. The generalized version, \cite{S}, allows for specialization of the form $x_i=(-1)^{s_i} a^{a_i} q^{q_i} x^{n_i}$, where $n_i$ are non-negative integers. In both versions the BPS numbers of the knot $K$ can be expressed as the integral linear combinations of the motivic Donaldson-Thomas of the corresponding quiver $Q$.

\section{Explicit formulae for the quiver generating series}\label{comb}

Let $Q$ be a symmetric quiver with $m$ vertices, and let $C$ be the corresponding adjacency matrix. Therefore, from now on $C$ is an integral symmetric $m\times m$ matrix. In this section we provide general expressions for coefficients $b_{l_1,\ldots,l_m}$ of the classical limit of the generating series 
\begin{equation}
y(x_1,\ldots,x_m)=\lim_{q\to 1} \frac{P_C(q^{n_1} x_1,\ldots,q^{n_m} x_m)}{P_C(x_1,\ldots,x_m)}=\sum_{l_1,\ldots,l_m} b_{l_1,\ldots,l_m} x_1^{l_1}\ldots x_m^{l_m},   \label{definv-bis}
\end{equation}
where $P_C(x_1,\ldots,x_m)$ is the motivic generating series introduced in (\ref{P-C}), and $n_i$ are arbitrary nonnegative integers.
These results naturally generalize the expressions obtained in \cite{PSS}, where the formulas were obtained in the case $n_1=\cdots=n_m=1$.

\begin{definition}\cite{PSS}
Let $k\in\{1,\ldots,m\}$. For a set
of $k$ pairs $(i_u,j_u)$, $u=1,\ldots,k$, where $1\le i_u,j_u\le m$, we say that it is admissible, if it satisfies the following two conditions:\\
(1)\quad  there are no two equal among $j_1,\ldots,j_k$\\
(2)\quad  there is no cycle of any length: for any $l$, $1\le l \le k$, there is no subset of $l$ pairs $(i_{u_{\ell}},j_{u_{\ell}})$, $\ell=1,\ldots,l$, such that $j_{u_{\ell}}=i_{u_{\ell+1}}$, $\ell=1,\ldots,l-1$, and $j_{u_{l}}=i_{u_1}$.
\end{definition}
Using this definition, we formulate the following proposition.
\begin{proposition}\label{glavna}
Constants $b_{l_1,\ldots,l_m}$ from (\ref{definv-bis}) take the form of
\begin{equation}
b_{l_1,\ldots,l_m} = (-1)^{\sum_{i=1}^m (C_{i,i}+1) l_i} A(l_1,\ldots,l_m)
\prod_{j=1}^m\frac{1}{n_j+\sum_{i=1}^m{C_{i,j}l_i}} \left({n_j+\sum_{i=1}^m{C_{i,j}l_i}} \atop {l_j}  \right),    \label{b1-glavna}
\end{equation}
where we define
\begin{equation}\label{forA1}
A(l_1,\ldots,l_m)= n_1\cdot\ldots \cdot n_m\left (1+\sum_{k=1}^{m-1}\sum_{admissible\, \Sigma_k}\prod_{(i_u,j_u)\in\Sigma_k} C_{i_u,j_u} l_{i_u}n_{j_u}^{-1}\right ).   
\end{equation}
Here, in the latter sum, we are summing over all admissible subsets of length $k$ -- one such subset we denote as $\Sigma_k$. 
\end{proposition}

Analogously to \cite{PSS}, Proposition \ref{glavna} can be proven by induction. As in \cite{PSS}, examples for small $m$ are presented in the following subsection, and (\ref{b1-glavna}) can be viewed as a generalization of those expressions. One can check that definition (\ref{forA1}) is formally valid even if some $n_i$ is zero; we treat $n_i$ as a nonzero integer and perform all the cancelations in the expression, substituting zero as a value of this variable at the end.

\subsection{Examples for $m=1,2$ and $3$}

Having obtained the general formula, we proceed to examples for small $m$. First, for $m=1$, we consider a quiver that consists of a single vertex with level $n_1=N$ and $f\in \Z$ loops, with a quiver matrix
\begin{equation}\label{1vert}
C = \left( f \right).
\end{equation}
Coefficients ${b}_{i}$ in (\ref{b1-glavna}) are then given by 
\begin{equation} 
{b}_{i} = \frac{(-1)^{(f+1) i}N}{f i +N}\left({f  i+N \atop i}\right).     \label{1vertcoeffb}
\end{equation}
Next, we consider a quiver with $m=2$ vertices with levels $n_1\equiv M$ and $n_2\equiv N$, with a $2\times 2$ symmetric quiver matrix
\begin{equation}\label{2vert}
C=\left(\begin{array}{cc}
\alpha & \beta \\
\beta & \gamma
\end{array}
\right),
\end{equation}
where $\alpha$, $\beta$ and $\gamma$ are arbitrary integers. Coefficients ${b}_{i,j}$ in (\ref{b1-glavna}) take the form of 
\begin{equation}
{b}_{i,j} = (-1)^{(\alpha+1) i + (\gamma+1) j}\frac{M\beta i+N\beta j+MN}{(\alpha i+\beta j+M)(\beta i+\gamma j+N)}\left({\alpha i+\beta j+M \atop i}\right)\left({\beta i+\gamma j+N \atop j}\right).   \label{C2x2-b-ij}
\end{equation}
Finally, we consider a quiver with $m=3$ vertices with level $n_1=M$, $n_2=N$ and $n_3=P$, whose $3\times 3$ symmetric quiver matrix is given by
\begin{equation}\label{trojka}
C=\left(\begin{array}{ccc}
\alpha & \beta &\delta \\
\beta & \gamma&\epsilon\\
\delta&\epsilon&\phi
\end{array}
\right),
\end{equation}
where $\alpha$, $\beta$, $\gamma$, $\delta,$ $\epsilon$, and $\phi$ are arbitrary integers. Coefficients ${b}_{i,j,k}$ in (\ref{b1-glavna}) are given by 
\begin{eqnarray*}
  {b}_{i,j,k} =&  (-1)^{(\alpha+1) i + (\gamma+1) j + (\phi+1) k} A_{i,j,k}\left({\alpha i+\beta j+\delta k+M \atop i}\right)\left({\beta i+\gamma j+ \epsilon k+N \atop j}\right)\nonumber \\
  &\times \left({\delta i+\epsilon j+ \phi k+P \atop k}\right),
\end{eqnarray*}
where we defined 
\begin{eqnarray*}
A_{i,j,k}=&\frac{1}{(\alpha i+\beta j+\delta k+M) (\beta i+\epsilon k+\gamma j+N) (\delta i+\epsilon j+\phi k+P)}\times\\&\left(\beta i j (N \delta+M\epsilon)+\delta i k (P \beta+M\epsilon)+\epsilon j k (P \beta+N \delta)+M\beta \delta i^2+i (P \beta+N \delta)\right.\\&+\left. N \beta \epsilon j^2+N j (P \beta+M\epsilon)+P \delta \epsilon k^2+P k (N \delta+M\epsilon)+MNP\right).
\end{eqnarray*}

We will use those expressions in the following sections. Note that Proposition \ref{glavna}, together with coefficients $b_{l_1,\dots,l_r}$ is, in a way, a minimal modification of the respective proposition and formulas for the level one generators.
\subsection{Relationship with the partial limits}

In \cite{PS}, in relationship with the computations of partition functions in strip geometries,  
the so-called partial limits were studied. 
\begin{equation}
y_j(x_1,\ldots,x_m)=\lim_{q\to 1} \frac{P_C(x_1,\ldots,x_{j-1},q x_j,x_{j+1},\ldots, x_m)}{P_C(x_1,\ldots,x_m)}=\sum_{l_1,\ldots,l_m} c^{(j)}_{l_1,\ldots,l_m} x_1^{l_1}\cdots x_m^{l_m},   \label{definv-bis23}
\end{equation}

The explicit form of the coefficients $ c^{(j)}_{l_1,\ldots,l_m} $ were obtained in \cite{PS}. By using these partial limits, we can express the limits of the quiver generating series with higher level generators in the following way
\begin{equation}
y(x_1,\ldots,x_m)=\lim_{q\to 1} \frac{P_C(q^{n_1}x_1,\ldots,q^{n_{j-1}}x_{j-1},q^{n_j} x_j,q^{n_{j+1}}x_{j+1},\ldots, q^{n_m}x_m)}{P_C(x_1,\ldots,x_m)}=   \label{definv-bis1}
\end{equation}
$$=y_1(x_1,\ldots,x_m)^{n_1}y_2(x_1,\ldots,x_m)^{n_2}\cdots y_m(x_1,\ldots,x_m)^{n_m}=\prod_{i=1}^m(y_i(x_1,\ldots,x_m))^{n_i}.$$

In such a way we can express the coefficients $b_{i_1,\ldots,i_m}$ of $y(x_1,\ldots,x_m)$ as a convolution of $ c^{(j)}_{l_1,\ldots,l_m}$. Together with our results from section \ref{comb}, this gives many combinatorial identities involving binomial coefficients, like  for example:
\begin{equation}
\frac{2}{f n +2}\left({f  n+2 \atop n}\right)=\sum_{i+j=n}\frac{1}{f i +1}\left({f  i+1 \atop i}\right)\frac{1}{f j +1}\left({f  j+1 \atop i}\right).
\end{equation}

\section{BPS-like numbers}

We have already introduced the BPS numbers in Section \ref{prelim}. From a physical perspective, one can use string theoretic arguments to show that they are integers. On the other hand, for any knot that has a respective quiver (in the sense of generalized knots-quivers correspondence), it is automatically guaranteed that BPS numbers are indeed integers \cite{S}. In this section, we devote our attention to calculating BPS numbers for those knots that have a quiver description using higher-level generators. Details of this calculation (for the case of level one generators) can be found in \cite{PSS}, but for completeness, we review the basic steps of the computation here. 

BPS numbers $N_{r,i,j}$ have three indices, where index $i$ is related to the $a$ variable, and index $j$ to the $q$ variable. We will be interested in the (semi)classical limit ($q\rightarrow 1$) and in the extremal limit (bottom row) in $a$ (though this will be problematic for those knots that do not obey the exponential growth property). Therefore, we will compute numbers $N_r$ that have only one index, and that are obtained by summing over the last index $j$ and appropriately identifying the second index with the multiple of the first one \cite{PSS}. From now on, we will call those numbers simply BPS numbers, or extremal BPS numbers. 

Of course, a general symmetric quiver is not guaranteed to have a knot representative, and therefore, we use the name BPS-like numbers when talking about integer invariants obtained from quivers without a specific knot in mind. To compute BPS numbers, we first compute the logarithm of $y(x)$ (in case of a quiver, we have $y(x_1,\dots,x_m)$ that gives $y(x)$ once appropriate specification of $x_i$ is made), and then the logarithmic derivative $ x(\log y(x))'$. Expanding this in power series in $x$, we have $x(\log y(x))'=\sum_i a_ix^i$. From coefficients $a_i$ we can compute BPS numbers as 
\begin{equation}\label{moe}
    N_r=\frac{1}{r^2}\sum_{d|r}\mu(d)a_{\frac{r}{d}},
\end{equation}
where $\mu(d)$ stands for the Moebious $\mu$ function.
First, let us compute the logarithm of $y(x_1,\dots ,x_m)$ from (\ref{definv-bis0}). Note that in the case without higher level generators, this was done in \cite{PSS}, where it was shown that the logarithm of $y$ takes the form $$\sum_{l_1,\ldots,l_m}c_{l_1,\ldots,l_m}x_1^{l_1}\cdots x_m^{l_m},$$ where coefficients $c_{l_1,\ldots,l_m}$ are closely related to previously introduced coefficients $b_{l_1,\ldots,l_m}$. It turns out that this observation stands also in the case of higher level generators. Namely, we have 

\begin{eqnarray}
     \log y(x_1,\dots ,x_m)&=&  \sum_{(l_1,\ldots,l_m)\neq {\bf 0}} (-1)^{\sum_{i=1}^m (C_{i,i}+1) l_i} A_{max}(l_1,\ldots,l_m) \times\label{logy}\\
&&\quad\quad\quad \times\prod_{j=1}^m  \frac{1}{\sum_{i=1}^m{C_{i,j}l_i}} \left({\sum_{i=1}^m{C_{i,j}l_i}} \atop {l_j}  \right)x_1^{l_1}\dots x_m^{l_m},\nonumber
\end{eqnarray}
where 
\begin{equation}
    A_{max}=A(l_1,\ldots,l_m)= n_1\cdot\ldots \cdot n_m\left (\sum_{admissible\, \Sigma_{m-1}}\prod_{(i_u,j_u)\in\Sigma_k} C_{i_u,j_u} l_{i_u}n_{j_u}^{-1}\right ). 
\end{equation}
The sum in (\ref{logy}) is taken over non-negative $l_i$'s, such that at least one of them is greater that zero. Note that a special care is required when computing the last expression, as for certain combinations of $l_i$ we may formally obtain expressions $\frac{0}{0}$,  and in that case we have to perform all necessary cancellations in the expression before computing with concrete values of $l_i$.

It is clear from the last equation that the only change between the last formula and the similar one from \cite{PSS} is the change of coefficients in the numerator of $A_{max}$. An important consequence of this is that the integer property of classical BPS numbers from \cite{PSS} implies the integrability of  BPS numbers in our case.

We can now compute BPS-like numbers (diagonal DT invariants), by specifying $x_i=x^{n_i}$. We have that the logarithmic derivative is given by 


\begin{eqnarray*}
     x(\log y(x))'&=&  \sum_{(l_1,\ldots,l_m)\neq {\bf 0}} (-1)^{\sum_{i=1}^m (C_{i,i}+1) l_i} A_{max}(l_1,\ldots,l_m) \times\\
&&\quad\quad\quad \times\,\,\prod_{j=1}^m  \frac{n_1l_1+\ldots+n_ml_m}{\sum_{i=1}^m{C_{i,j}l_i}} \left({\sum_{i=1}^m{C_{i,j}l_i}} \atop {l_j}  \right)x^{n_1l_1+\ldots+ n_m l_m}.
\end{eqnarray*}

We can compute the BPS like numbers from the last formula using a well-known procedure, as in \cite{KRSS}, i.e. by (\ref{moe}). We will do this explicitly in the next subsection, using concrete examples of knots. 

For concreteness, let us explicitly write the case of the single vertex quiver (\ref{1vert}) of level $N$, with $f$ loops. The logarithmic derivative is then given by 
\begin{equation}\label{jedan}
     x(\log y(x))'=\sum_{l}(-1)^{(f+1)l}\frac{N^2}{f}\binom{fl}{l}x^{Nl}.
\end{equation}
As already anticipated, the integrability of BPS numbers is guaranteed from the integrability of BPS numbers for a level one quiver, due to the fact that $N$ is an integer. We will use this expression in the next subsection, where we will provide some concrete examples.

\subsection{Knots examples}


As explained in the introduction, the main motivation for our work comes from the fact that certain knots, upon using knots-quivers correspondence, have a quiver generating series with higher level generators \cite{S}. Let us not work out an example of $\mathbf{9}_{46}$ knot, whose reduced $S^r$ colored HOMFLY-PT polynomial is given \cite{S} . Focusing on a bottom row (lowest order in $a$), we have 
\begin{equation}
P_r(q)=\sum_{d_1+d_3+2d_4=r}\frac{(q;q)_r q^{-d_1 (d_3+2d_4)-2 d_4 (d_3+d_4)}}{(q;q)_{d_1}
   (q;q)_{d_3} (q;q)_{d_4}}.
\end{equation}
For example, we have $P_1(q)=2$, $P_2(q)=2+\frac{1}{q}+\frac{1}{q^2}$, $P_3(q)=\frac{2}{q^4}\left(1+q+q^2+q^4\right)$, and so on. Generating function is therefore given by 
\begin{equation}
    P(x,q)=\sum_{d_1,d_3,d_4} \frac{(q;q)_{d_1+d_3+2d_4} q^{-d_1 (d_3+2d_4)-2 d_4 (d_3+d_4)}}{(q;q)_{d_1}
   (q;q)_{d_3} (q;q)_{d_4}} x^{d_1} x^{d_2}x^{2d_4}.
\end{equation}
Now we make the transition to the unreduced version of HOMFLY-PT polynomials so that we get 
\begin{equation}
     \overline{P}(x,q)=\sum_{d_1,d_3,d_4} \frac{ q^{-d_1 (d_3+2 d_4)+\frac{1}{2} (d_1+d_3+2 d_4)-2 d_4
   (d_3+d_4)}}{(q;q)_{d_1}
   (q;q)_{d_3} (q;q)_{d_4}} x^{d_1} x^{d_2}x^{2d_4}.
\end{equation}
From this expression, we obtain that symmetric $3\times 3$ matrix $C$ is given by 
\begin{equation}
   C= \begin{pmatrix}
0&-1&-2\\
-1&0&-2\\
-2&-2&-4\\
    \end{pmatrix},
\end{equation}
and it is obvious that we have one generator of level two ($M=1, N=1$ and $P=2$ in (\ref{trojka})).
Logharitmic derivative is 
\begin{equation}
  x(\log y(x))'=  -2 x - 10 x^2 - 56 x^3 - 330 x^4+\dots
\end{equation}
For this particular case, we can compute BPS numbers. The first few are given by $N_1=-2$, $N_2=-2$, $N_3=-6$, and $N_4=-20$, and they are obviously integers, as expected.

We now move to $\mathbf{8}_{20}$ knot, using again the result for the unreduced HOMFLY-PT polynomial from \cite{S}. Focusing again on the bottom row, we have that the $\overline{P}_r(q)$ is equal to
\begin{equation}
   \sum_{d_1+d_2+d_3+2d_4=r} \frac{(-1)^{d_2+d_3+d_4}q^{\frac{1}{2}(-d_2^2+d_3^2+d_3(1-2d_1-2d_4)-d_4(1+4d_1+5d_4)-d_2(1+2d_1+2d_3+6d_4))+\frac{r}{2}}}{(q;q)_{d_1}(q;q)_{d_2}(q;q)_{d_3}(q;q)_{d_4}}.
\end{equation}
Similar to the previous example, we have one generator of level two, but here we have three generators of level one. The respective symmetric quiver matrix is 
\begin{equation}
    C=\begin{pmatrix}
        0&-1&-1&-2\\
        -1&-1&-1&-3\\
        -1&-1&1&-1\\
        -2&-3&-1&-5\\
    \end{pmatrix}.
\end{equation}
The explicit expression of the logarithmic derivative is given by  
\begin{equation}
    x(\log y(x))'=x+5x^2-17x^3+5x^4+\dots,
\end{equation}
from which we calculate BPS numbers as $N_1=1$, $N_2=1$, $N_3=-2$, $N_4=0$... As expected, they are integers.

 Another well-known example is knot $\mathbf{9}_{42}$. This knot is notorious for not obeying the exponential growth property. In particular, the non-zero terms of the lowest degree in $a$ of the colored HOMFLY-PT polynomial $P_r(a,q)$ is not linear function in $r$. Rather, we can use the following extremal specialization 
 $$P_r(q)=a^{\frac{3}{2}r}P_r(a,q)_{{\mid}_{a\to 0}}.$$ 
 In such a way, we get the extremal (bottom-row like) specializations of the colored HOMFLY-PT polynomial for the knot $\mathbf{9}_{42}$, such that $P_r(q)$ is equal to zero for odd $r$, and nonzero for even $r$.
 Therefore, we can focus on even degrees of $r=2l$, where the unreduced version of the polynomial is given by \cite{S}
\begin{equation}
    \overline{P}(x,q)=\sum_l q^{-6l^2+2l}\frac{q^l}{(q;q)_l}x^{2l} ,
\end{equation}
so that the quiver matrix we get is 
\begin{equation}
    C=(-12),
\end{equation}
with $n_1\equiv N=2$.
Using relation (\ref{jedan}), we have 
\begin{equation}
     x(\log y(x))'=-4x^2-100x^4-2812x^6-83300x^8+\dots.
\end{equation}

Classical BPS-like numbers are given by $N_1=0$, $N_2=-1$, $N_3=0$, $N_4=-6$, $N_5=0$, $N_6=-78$, $N_7=0$, $N_8=-1300$,... Again, it is clear that those numbers are integers. 
\section{Combinatorial interpretation}\label{komb}
\subsection{Lattice paths under lines of integer slope}

One of the numerous famous interpretations of the Catalan numbers is the count of the lattice paths in the positive quadrant in the plane, below the line $y=x$. A natural generalization of this is provided by Fuss-Catalan numbers, which count the lattice paths below the line $y=ax$, for some non-negative integer $a$. Fuss-Catalan numbers where computed in \cite{PSS} by using quiver generating series of a single vertex quiver with $a+1$ loops.

The single vertex quiver generating series, when the vertex is of higher level, turns out to be related to lattice paths under the line $y=ax+b$, where both $a$ and $b$ are non-negative integers. 
In such a way we obtain the count of lattice paths under the critical line that doesn't necessarily pass through the origin, which is one of the novelties in the relationship between quivers and enumerative combinatorics of lattice paths.

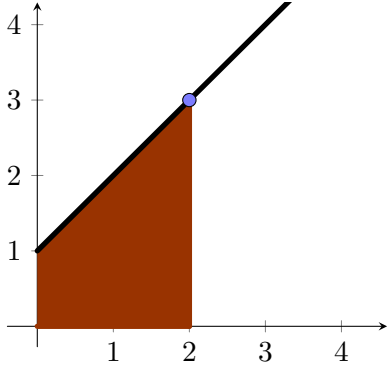
\begin{figure}[h!]
    \centering
   \begin{tikzpicture}[line cap=round,line join=round,>=triangle 45,x=1.0cm,y=1.0cm]
\begin{axis}[
x=1.0cm,y=1.0cm,
axis lines=middle,
xmin=-0.393641438992434,
xmax=4.60077465954524,
ymin=-0.270664752876202,
ymax=4.2952857862405,
xtick={-6.0,-5.0,...,7.0},
ytick={-3.0,-2.0,...,6.0},]
\clip(-6.393641438992434,-3.270664752876202) rectangle (7.60077465954524,6.2952857862405);
\fill[line width=2.pt,color=zzttqq,fill=zzttqq,fill opacity=0.10000000149011612] (0.,0.) -- (2.,0.) -- (2.,3.) -- (0.,1.) -- cycle;
\draw [line width=2.pt,domain=0.0:7.60077465954524] plot(\x,{(--1.--1.*\x)/1.});
\draw [line width=2.pt,color=zzttqq] (0.,0.)-- (2.,0.);
\draw [line width=2.pt,color=zzttqq] (2.,0.)-- (2.,3.);
\begin{scriptsize}
\draw [fill=xdxdff] (2.,3.) circle (2.5pt);
\end{scriptsize}
\end{axis}
\end{tikzpicture}
\caption{Example of a path in the case $a=1$, $b=1$, where  the area of the shaded region is equal to four. }
\label{fig1}
\end{figure}

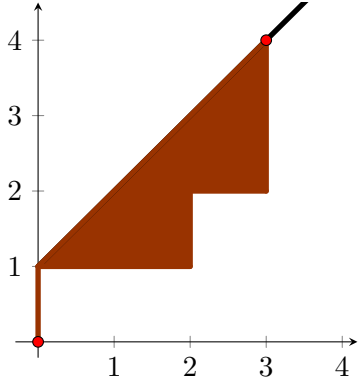
\begin{figure}[h!]
\centering
  \begin{tikzpicture}[line cap=round,line join=round,x=1.0cm,y=1.0cm]
\begin{axis}[
x=1.0cm,y=1.0cm,
axis lines=middle,
xmin=-0.294081433485506,
xmax=4.2,
ymin=-0.205743226091135,
ymax=4.5,
xtick={-5.0,-4.0,...,8.0},
ytick={-4.0,-3.0,...,5.0},]
\clip(-5.294081433485506,-4.205743226091135) rectangle (8.873910852009836,5.478856151784831);
\fill[line width=2.pt,color=zzttqq,fill=zzttqq,fill opacity=0.10000000149011612] (0.,1.) -- (2.,1.) -- (2.,2.) -- (3.,2.) -- (3.,4.) -- cycle;
\draw [line width=2.pt,domain=0.0:8.873910852009836] plot(\x,{(--2.--2.*\x)/2.});
\draw [line width=2.pt,color=zzttqq] (0.,1.)-- (2.,1.);
\draw [line width=2.pt,color=zzttqq] (2.,1.)-- (2.,2.);
\draw [line width=2.pt,color=zzttqq] (2.,2.)-- (3.,2.);
\draw [line width=2.pt,color=zzttqq] (3.,2.)-- (3.,4.);
\draw [line width=2.pt,color=zzttqq] (3.,4.)-- (0.,1.);
\draw [line width=2.pt,color=zzttqq] (0.,0.)-- (0.,1.);
\begin{scriptsize}
\draw [fill=ffqqqq] (0.,0.) circle (2.0pt);
\draw [fill=ffqqqq] (3.,4.) circle (2.0pt);
\end{scriptsize}
\end{axis}
\end{tikzpicture}
\caption{Example of a path between $(0,0)$ and $(3,4)$, case of $a=1$ and $b=1$. The shaded area is the area of the path, and in a given case, equals to $\frac{7}{2}$.}\label{fig2}
\end{figure}

To that end, let $a$ and $b$ be nonnegative integers. Let $\mathcal{S}_n$ be the set of all  lattice paths from $(0,0)$ to $(n,an+b)$, that stay below the critical line $y=ax+b$, and let $A_n$ be the number of such paths. By $A_n(q)$ we denote the weighted count of such lattice paths, that is, to each such path $\tau$ we associate weight $q^{\mathrm{Area}(\tau)}$, where $\mathrm{Area}(\tau)$ denotes the area between the critical line and the path $\tau$ (see examples on Figure \ref{fig1} and Figure \ref{fig2}). We will first work out the full quantum case and later specialize to the classical $q\rightarrow 1$ limit.

Now, every path $\tau$ is determined bijectively by the $y$-coordinates of the horizontal parts of the lattice paths. Denote by $y_i$, $i=1,\ldots,n$, the $y$-coordinates of horizontal part of $\tau$ for $x$-coordinate satisfying $i-1\le x\le i$. In other words, path $\tau$ is of the form $(0,0)$--$(0,y_1)$--$(1,y_1)$--$(1,y_2)$--$(2,y_2)$--$(2,y_3)$--$\ldots$--$(n-1,y_n)$--$(n,y_n)$--$(n,ax+b)$, where coordinates $y_i\ge 0$, $i=1,\ldots,n$, satisfy $y_i\le a(i-1)+b$, and $y_i\le y_{i+1}$, for $i=1,\ldots,n-1$. Equivalently this can be rewritten in terms of variables $c_i:=a(i-1)+b-y_i$, $i=1,\ldots,n$, that measure the distance from the critical line.

In such a way, every lattice path $\tau\in\mathcal{S}_n$ is uniquely determined by the sequence $(c_1,\ldots,c_n)$ of nonnegative integers, such that $c_1\le b$, and $c_{i+1}\le c_i+a$, for $i=1,\ldots,n-1$.
Therefore 
\begin{equation}\label{An}
A_n=\sharp\{(c_1,\ldots,c_n)\mid 0\le c_1\le b,\, 0\le c_{i+1}\le c_i+a,\,i=1,\ldots,n-1\}.
\end{equation}
Also, the area of the path $\tau$ is then equal to $\frac{an}{2}+\sum_{i=1}^n c_i$, and so
\begin{equation}\label{Aqn}
A_n(q)=q^{\frac{1}{2}a n}\sum_{c_1=0}^{b}\sum_{c_2=0}^{c_1+a}\sum_{c_3=0}^{c_2+a}\cdots\sum_{c_n=0}^{c_{n-1}+a} {q^{c_1+c_2+c_3+\ldots+c_{n}}}.
\end{equation}

\subsection{Single level $N$ generator}
The numbers $A_n(q)$, and consequently $A_n$, naturally appear in the growth rate quotient of the quiver generating series for the single vertex quiver with higher level generator.

\begin{proposition}\label{levN}
Let $f\ge 0$ be a nonnegative integer, and let $N\ge 1$ be a positive integer. Let
\begin{equation}
P_N(q,x)=\sum_{d\ge 0} (-1)^{f d}\frac{q^{\frac{1}{2}f(d^2-d)}}{(q;q)_d} x^{Nd},
\end{equation}
be a quiver generating series corresponding to a quiver with a single generator of level $N$, and with $f$ self-loops. Then
\begin{equation}\label{f12}
\frac{P_N(q,qx)}{P_N(q,x)}=1+\sum_{n>0} a^{(N)}_n(q) x^{Nn},
\end{equation}
where the coefficients $a^{(N)}_n(q)$ are given by:
\begin{equation}\label{foran}
a^{(N)}_n(q)=(-1)^{(f-1)n}\sum_{c_1=0}^{N-1}\sum_{c_2=0}^{c_1+f-1}\sum_{c_3=0}^{c_2+f-1}\sum_{c_4=0}^{c_3+f-1}\cdots\sum_{c_n=0}^{c_{n-1}+f-1} {q^{c_1+c_2+c_3+c_4+\ldots+c_{n}}}, \quad n\ge 1.
\end{equation}
\end{proposition} 

Before proving the proposition, we give a simple lemma that enables rewriting the general term of the generating series $P_N(q,x)$ as a suitable multiple sum
\begin{lemma}\label{lem12}
For every $d\ge 1$, and $|q|<1$, one has:
\begin{equation}\label{flem12}
\frac{q^{\frac{1}{2}f(d^2-d)}}{(q;q)_d}=\sum_{c_1\ge 0}\,\sum_{c_2\ge c_1+f}\,\sum_{c_3\ge c_2+f}\cdots\sum_{c_d\ge c_{d-1}+f}\, {q^{c_1+c_2+c_3+\ldots+c_{d}}}.
\end{equation}
\end{lemma}

\textbf{Proof:} Expression (\ref{flem12}) follows by induction on $d$. Indeed, for $d=1$ it is just
$$\frac{1}{1-q}=\sum_{c_1\ge 0}\,q^{c_1}.$$
As for that induction step, supposing that (\ref{flem12}) holds for some $d$, we have
\begin{eqnarray*}
\!\!\!\!\!\!\!\frac{q^{\frac{1}{2}f((d+1)^2-(d+1))}}{(q;q)_{d+1}}&=&\frac{q^{fd}}{1-q^{d+1}}\,\frac{q^{\frac{1}{2}f(d^2-d)}}{(q;q)_d}=\\
&=&q^{fd}\sum_{e\ge 0} q^{(d+1)e} \sum_{c_1\ge 0}\,\sum_{c_2\ge c_1+f}\,\sum_{c_3\ge c_2+f}\cdots\sum_{c_d\ge c_{d-1}+f}\, {q^{c_1+c_2+c_3+\ldots+c_{d}}}=\\
&=&\sum_{e\ge 0}\sum_{c_1\ge 0}\,\sum_{c_2\ge c_1+f}\,\cdots\sum_{c_d\ge c_{d-1}+f}\, {q^{e+(c_1+e+f)+(c_2+e+f)\ldots+(c_{d}+e+f)}}=\\
&=&\sum_{c'_1\ge 0}\,\sum_{c'_2\ge c'_1+f}\,\sum_{c'_3\ge c'_2+f}\,\cdots\sum_{c'_{d+1}\ge c'_{d}+f}\, {q^{c'_1+c'_2+c'_3+\ldots+c'_{d+1}}},
\end{eqnarray*}
where we set $c'_1=e$, and $c'_{i+1}=c_i+e+f$, for $i=1,\ldots,d$. \kraj

\textbf{Proof of Proposition \ref{levN}:} Let $b_d=(-1)^{f d}\frac{q^{\frac{1}{2}f(d^2-d)}}{(q;q)_d}$, $d\ge 0$, and so $P_N(x)=\sum_{d\ge 0} b_d x^{Nd}$. Then the statement of the Proposition is equivalent to:
\begin{equation}\label{for13}
q^{n N} b_n=a^{(N)}_n+b_1 a^{(N)}_{n-1}+b_2 a^{(N)}_{n-2}+\ldots+b_n.
\end{equation} By using Lemma \ref{lem12} to express $b_d$, as well as (\ref{foran}), the right-hand side of (\ref{for13}) becomes:
\begin{eqnarray*}
\!\!\!(-1)^{(f-1) n}\!\!\!&\!\!&\!\!\!\left(\sum_{c_1=0}^{N-1}\sum_{c_2=0}^{c_1+f-1}\cdots\sum_{c_{n-1}=0}^{c_{n-2}+f-1}\sum_{c_n=0}^{c_{n-1}+f-1} {q^{\sum_{i=1}^n c_i}}-\sum_{c_1=0}^{N-1}\sum_{c_2=0}^{c_1+f-1}\cdots\sum_{c_{n-1}=0}^{c_{n-2}+f-1}\sum_{c_n\ge 0} {q^{\sum_{i=1}^n c_i}}+\right.\\
&&\quad+\sum_{c_1=0}^{N-1}\sum_{c_2=0}^{c_1+f-1}\cdots\sum_{c_{n-2}=0}^{c_{n-3}+f-1}\sum_{c_{n-1}\ge 0}\sum_{c_n\ge c_{n-1}+f} {q^{\sum_{i=1}^n c_i}}-\ldots\\
&&\left.\quad\ldots+(-1)^{n}\sum_{c_1\ge 0}\,\sum_{c_2\ge c_1+f}\cdots\sum_{c_{n-1}\ge c_{n-2}+f}\,\,\sum_{c_n\ge c_{n-1}+f} {q^{\sum_{i=1}^n c_i}}\right)=\end{eqnarray*}
\begin{eqnarray*}
&&=(-1)^{(f-1) n}\left(-\sum_{c_1=0}^{N-1}\sum_{c_2=0}^{c_1+f-1}\cdots\sum_{c_{n-1}=0}^{c_{n-2}+f-1}\sum_{c_n\ge c_{n-1}+f} {q^{\sum_{i=1}^n c_i}}+\right.\\
&&\quad+\sum_{c_1=0}^{N-1}\sum_{c_2=0}^{c_1+f-1}\cdots\sum_{c_{n-2}=0}^{c_{n-3}+f-1}\sum_{c_{n-1}\ge 0}\sum_{c_n\ge c_{n-1}+f} {q^{\sum_{i=1}^n c_i}}-\ldots\\
&&\left.\quad\ldots+(-1)^{n}\sum_{c_1\ge 0}\,\sum_{c_2\ge c_1+f}\cdots\sum_{c_{n-1}\ge c_{n-2}+f}\,\,\sum_{c_n\ge c_{n-1}+f} {q^{\sum_{i=1}^n c_i}}\right)=\end{eqnarray*}
\begin{eqnarray*}
&&=(-1)^{(f-1) n}\left(-\sum_{c_1=0}^{N-1}\sum_{c_2=0}^{c_1+f-1}\cdots\sum_{c_{n-2}=0}^{c_{n-3}+f-1}\sum_{c_{n-1}\ge c_{n-2}+f}\sum_{c_n\ge c_{n-1}+f} {q^{\sum_{i=1}^n c_i}}+\ldots\right.\\
&&\left.\quad\ldots+(-1)^{n}\sum_{c_1\ge 0}\,\sum_{c_2\ge c_1+f}\cdots\sum_{c_{n-1}\ge c_{n-2}+f}\,\,\sum_{c_n\ge c_{n-1}+f} {q^{\sum_{i=1}^n c_i}}\right)=\ldots\\
&&\ldots = (-1)^{fn} \sum_{c_1\ge N}\,\sum_{c_2\ge c_1+f}\,\sum_{c_3\ge c_2+f}\cdots\sum_{c_n\ge c_{n-1}+f}\, {q^{\sum_{i=1}^n c_i}}.
\end{eqnarray*}
Finally the last expression is equal to $q^{nN} b_n$ since
\begin{eqnarray*}
q^{nN} b_n&=&(-1)^{fn}q^{nN}\sum_{c_1\ge 0}\,\sum_{c_2\ge c_1+f}\,\sum_{c_3\ge c_2+f}\cdots\sum_{c_n\ge c_{n-1}+f}\, {q^{\sum_{i=1}^n c_i}}=\\
&=&(-1)^{fn} \sum_{c_1\ge 0}\,\sum_{c_2\ge c_1+f}\,\sum_{c_3\ge c_2+f}\cdots\sum_{c_n\ge c_{n-1}+f}\, {q^{\sum_{i=1}^n (c_i+N)}}=\\
&=&(-1)^{fn} \sum_{c'_1\ge N}\,\sum_{c'_2\ge c'_1+f}\,\sum_{c'_3\ge c'_2+f}\cdots\sum_{c'_n\ge c'_{n-1}+f}\, {q^{\sum_{i=1}^n c'_i}},
\end{eqnarray*}
which finishes the proof of (\ref{for13}), as desired.
 \kraj

\subsection{Specialization $q\to 1$}

The above quotients in the classical $q\to 1$ limit give rise to algebraic functions:
$$y_1(x)=\frac{P_1(q,qx)}{P_1(q,x)}=1+\sum_{n\ge 1}a_nx^N$$ and 
$$y_N(x)= \frac{P_N(q,qx)}{P_N(q,x)}=1+\sum_{n\ge 1}a_n^{(N)}x^{Nn}$$

Here $$a_n=\sharp\{(c_1,c_2,\ldots,c_n)|c_1=0, 0\le c_{i+1}\le c_{i}+f-1\},$$
and $$a_n^{(N)}=\sharp\{(c_1,c_2,\ldots,c_n)|0\le c_1\le N-1, 0\le c_{i+1}\le c_{i}+f-1\}.$$

In particular, for $f=2$, the above expression for $a_n$ gives one of the combinatorial interpretations of Catalan numbers. For arbitrary $f$ these are Fuss-Catalan numbers, i.e. $a_n=\frac{1}{fn+1}\binom{fn+1}{n}$. For arbitrary $N$ we have $a_n^{(N)}=\frac{N}{fn+N}\binom{fn+N}{n}$. 
More directly, since $P_N(x)=P_1(x^N)$, we have that $$y_N(x)=y_1(x^N)^N.$$

\begin{lemma}
For arbitrary positive integers $N,M,c$ and nonnegative integer $m$, the following identity holds:
\begin{equation}
\sum_{j+k=m} \frac{N}{cj+N} \binom{cj+N}{j} \frac{M}{ck+M} \binom{ck+M}{k}=\frac{N+M}{cm+N+M}\binom{cm+N+M}{m}.
\end{equation}
\end{lemma}

\textbf{Proof:} Let $y=y(x)$ be the power series given by:
$$y(x)=\sum_{j\ge 0} \frac{1}{cj+1} \binom{cj+1}{j} x^j,$$
i.e. the algebraic function given by 
$$ xy(x)^c-y(x)+1=0,$$
and $y(0)=1$. Then
$$y(x)^N=\sum_{j\ge 0} \frac{N}{cj+N} \binom{cj+N}{j} x^j,$$ 
and the wanted equality follows from expanding the identity
$$y(x)^{N+M}=y(x)^N y(x)^M,$$
as a power series in two different ways and comparing the coefficients. \kraj

Note that the above Lemma also follows from the rewriting our classical limit as product of partial limits (\ref{definv-bis1}).

\subsection{Non-weighted count and Raney numbers}
\begin{figure}[h!]
    \centering

\begin{tikzpicture}[line cap=round,line join=round,x=1.0cm,y=1.0cm, scale=1.5]
\begin{axis}[
x=1.0cm,y=1.0cm,
axis lines=middle,
xmin=-0.05,
xmax=3.5,
ymin=-0.05,
ymax=3.5,
xtick={0,1,...,3},
ytick={0,1,...,3}]
\clip(-4.1094878456836685,-3.7248353433427246) rectangle (7.245656730451171,5.037028432536828);
\draw [line width=1.pt,domain=0.0:7.245656730451171] plot(\x,{(--1.--1.*\x)/1.});
\draw [line width=1.pt,color=zzttqq] (0.,0.)-- (2.,0.);
\draw [line width=1.pt,color=zzttqq] (2.,0.)-- (2.,3.);
\draw [line width=1.pt,color=ffqqtt] (0.,0.)-- (1.,0.);
\draw [line width=1.pt,color=ffqqtt] (1.,0.)-- (1.,1.);
\draw [line width=1.pt,color=ffqqtt] (1.,1.)-- (2.,1.);
\draw [line width=1.pt,color=wwccff] (0.,0.)-- (1.,0.);
\draw [line width=1.pt,color=wwccff] (1.,0.)-- (1.,2.);
\draw [line width=1.pt,color=wwccff] (1.,2.)-- (2.,2.);
\draw [line width=1.pt,color=wwccff] (2.,2.)-- (2.,3.);
\draw [line width=1.pt,color=qqwuqq] (0.,1.)-- (2.,1.);
\draw [line width=1.pt,color=qqwuqq] (2.,1.)-- (2.,3.);
\draw [line width=1.pt,color=qqwuqq] (0.,1.)-- (0.,0.);
\draw [line width=1.pt,color=xdxdff] (0.,1.)-- (1.,1.);
\draw [line width=1.pt,color=xdxdff] (1.,1.)-- (1.,2.);
\draw [line width=1.pt,color=xdxdff] (1.,2.)-- (2.,2.);
\draw [line width=1.pt,color=xdxdff] (2.,2.)-- (2.,3.);
\draw [line width=1.pt,color=xdxdff] (0.,1.)-- (0.,0.);
\begin{scriptsize}
\draw [fill=red] (2.,3.) circle (1.5pt);
\end{scriptsize}
\end{axis}

\end{tikzpicture}
\caption{Case $N=1$, $f=1$; different paths from $(0,0)$ to $(2,3)$. The total number of paths can be read from the picture and is equal to five, as easily confirmed by the formula (\ref{kombre}). }
\label{fig3}
\end{figure}
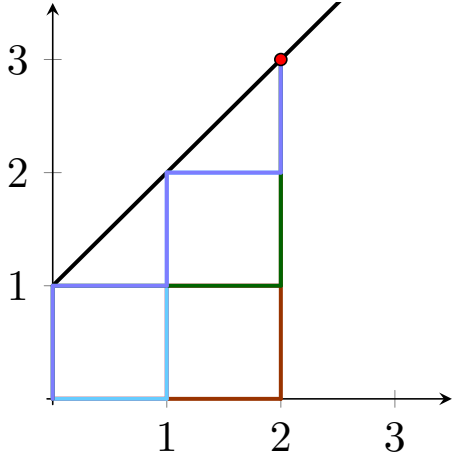

Let $C$ be a single vertex quiver with $f+1$ loops, with $f\ge 1$. 
Also, suppose the level of a generator is $N+1$, 
Then, as we have computed above, the coefficients of the corresponding classical limit are given by:
\begin{equation} 
{b}_{i} = (-1)^{(f+2) i}\frac{N+1}{(f+1)i +N+1}\left({(f+1)  i+N+1 \atop i}\right).  
\end{equation}
Up to an overall sign factor, these numbers are equal to
\begin{equation}\label{kombre}
     \frac{N+1}{(f+1)i +N+1}\left({(f+1)  i+N+1 \atop i}\right),
\end{equation}
which are the so-called Raney numbers. These numbers count precisely the number of lattice path in the positive quadrant consisting of right steps (1,0), and up steps (0,1), that start at (0,0) and end at $(i,fi+N)$, and that never go above the line $ y = fx+N$. Example is presented in Figure \ref{fig3}, and one can easily count the total number of paths to be five, as directly calculated from the formula (\ref{kombre}).\\

\section{Conclusion}
We present the asymptotics underlying the quiver generating series where generators can have higher levels, and extend the previous results when all generators where of level one. We also used this to compute the corresponding extremal BPS numbers for some of the thick knots with super-exponential growth, property which were inaccessible by previous methods.

On the combinatorial side, the coefficients in the corresponding generating series have interpretation as the count of lattice paths under the critical line. Here, for the first time the critical line does not need to pass through the origin, and the precisely the level $N$, measures the ``offset" of the critical line and its distance from the origin.

\section*{Acknowledgements}

Work of D.D. is supported by the funding provided by the Faculty of Physics, University of Belgrade, through grant number 451-03-47/2023-01/200162 by the Ministry of Science, Technological Development and Innovations of the Republic of Serbia and by the Science Fund of the Republic of Serbia, grant number TF C1389-YF, Towards a Holographic Description of Noncommutative Spacetime: Insights from Chern-Simons Gravity, Black Holes and Quantum Information Theory - HINT.  Work of M.S. has been supported by Funda\c{c}\~ao para a Ci\^encia e a Tecnologia (FCT) through
CEEC grant with DOI no. 10.54499/2020.02453.CEECIND/CP1587/CT0007, and by the Science Fund of the Republic of Serbia, Project no. 7749891, GWORDS -- ``Graphical Languages".



\end{document}